\documentclass[12pt]{amsart}

\usepackage{amsmath}
\usepackage{amsthm}
\usepackage{amsfonts}
\usepackage{amssymb}
\usepackage{amscd}

\usepackage[dvips]{graphicx}
\usepackage{psfrag}

\input{psfig}

\theoremstyle{plain}
\newtheorem{mainthm}{Theorem}

\newtheorem{thm}{Theorem}[section]

\newtheorem{lemma}[thm]{Lemma}

\newtheorem{defi}[thm]{Definition}

\newtheorem{rk}[thm]{Remark}
\newtheorem{q}[thm]{Question}

\title{Omega-limit sets
close to singular-hyperbolic attractors}

\author{C. M. Carballo and C. A. Morales}
\thanks{
The second author was partially supported by CNPq, FAPERJ and PRONEX-Dyn. Sys./Brazil.}

\begin{document}

\maketitle

\begin{abstract}
We study the omega-limit sets $\omega_X(x)$ in an isolating block $U$
of a singular-hyperbolic attractor
for three-dimensional vector fields $X$.
We prove that for every
vector field $Y$ close to $X$
the set $
\{x\in U:\omega_Y(x)$ contains a singularity$\}$
is {\em residual} in $U$.
This is used to prove the persistence of
singular-hyperbolic attractors
with only one singularity
as chain-transitive Lyapunov stable sets.
These results generalize well known properties of the
geometric Lorenz attractor \cite{gw}
and the example in \cite{mpu}.
\end{abstract}

\section{Introduction}

The {\em omega-limit set} of
$x$ with respect to
a vector field $X$ with generating flow $X_t$
is the accumulation point set
$\omega_X(x)$ of the positive orbit of $x$, namely
$$
\omega_X(x)=\left\{y:\,\,y=\lim_{t_n\to\infty}X_{t_n}(x)
\mbox{ for some sequence } t_n\to\infty\right\}.
$$
The structure of the omega limit sets
is well understood for vector fields
on compact surfaces.
In fact, the {\em Poincar\'e-Bendixon Theorem}
asserts that the omega-limit set
for vector fields with finite many singularities
in $S^2$ is either
a periodic orbit or a singularity
or a graph (a finite union of singularities
an separatrices forming a closed curve).
The {\em Schwartz Theorem}
implies that the omega-limit set
of a $C^\infty$ vector field
on a compact surface either contains
a singularity or an open set or is
a periodic orbit.
Another result is the {\em Peixoto Theorem}
asserting that
an open dense subset of vector fields
on any closed orientable surface are
{\em Morse-Smale},
namely their nonwandering set
is formed by a finite union of closed orbits
all of whose invariant manifolds are in general position.
A direct consequence this result
is that, for an open-dense subset of
vector fields on closed orientable surfaces,
most omega-limit sets are contained
in the attracting closed orbits.
This provides a complete description of
the omega limit sets on closed orientable surfaces.

The above results are known to be false
in dimension $>2$.
Hence
extra hypotheses
to understand the omega-limit sets are needed in general.
An important one is the hyperbolicity
introduced by Smale in the sixties.
Recall that a compact invariant set
is {\em hyperbolic} if it exhibits
contracting and expanding direction
which together with the flow's direction
form a continuous tangent bundle decomposition.
This definition leads the concept of
{\em Axiom A vector field},
namely the ones whose non-wandering set is
both hyperbolic and the closure of its closed orbits.
The Spectral Decomposition Theorem describes
the non-wandering set for
Axiom A vector fields, namely it decomposes
into a finite disjoint union of hyperbolic basic sets.
A direct consequence of the Spectral Theorem
is that for every Axiom A vector field
$X$ there is an open-dense subset
of points whose omega-limit set
are contained in the hyperbolic attractors
of $X$.
By {\em attractor} we mean
a compact invariant set
$\Lambda$ which is {\em transitive} (i.e. $\Lambda=\omega_X(x)$ for some $x\in \Lambda$)
and satisfies $\Lambda=\cap_{t\geq 0}X_t(U)$ for some
compact neighborhood $U$ of it called {\em isolating block}.
On the other hand,
the structure of the omega-limit sets
in an isolating block $U$ of a hyperbolic attractor
is well known: For every vector field
$Y$ close to $X$
the set
$$
\{x\in U:
\omega_Y(x)=\cap_{t\geq 0}Y_t(U)\}
$$
is {\em residual} in $U$.
In other words, the omega-limit sets
in a residual subset of $U$
are uniformly distributed in the maximal invariant set
of $Y$ in $U$.
This result is a direct consequence
of the structural stability of the
hyperbolic attractors.

There are many examples
of non-hyperbolic vector fields
$X$ with a large set
of trajectories going to the attractors of $X$.
Actually, a conjecture by Palis \cite{p}
claims that this is true for a dense set
of vector fields on any compact manifold
(although he used a different definition of attractor).
A strong evidence
is the fact that
there is a residual subset of $C^1$ vector fields
$X$ on any compact manifold exhibiting a residual subset
of points whose omega-limit sets
are contained in the chain-transitive Lyapunov stable sets of $X$ (\cite{mpa4}).
We recall that a compact invariant set $\Lambda$ is
{\em chain-transitive}
if any pair of points on it can be
joined by a pseudo-orbit with arbitrarily small
jump. In addition, $\Lambda$ is {\em Lyapunov stable} if
the positive orbit of a point close to $\Lambda$
remains close to $\Lambda$.
The result \cite{mpa4} is weaker than the Palis
conjecture since
every attractor is a chain-transitive Lyapunov stable set but not vice versa. 

In this paper
we study the omega-limit sets in an isolating block
of an attractor for vector fields on compact three
manifolds.
Instead of hyperbolicity we shall assume that
the attractor is
{\em singular-hyperbolic}, namely it
has singularities (all hyperbolic)
and is partially hyperbolic with volume expanding central direction \cite{mpp1}.
These attractors were considered in \cite{mpp1}
for a characterization of $C^1$ robust transitive sets with singularities for vector fields on compact three manifolds (see also \cite{mpp3}).
The singular-hyperbolic attractors are not hyperbolic
although they have some properties resembling the hyperbolic ones.
In particular, they
do not have the pseudo-orbit tracing property
and are neither
expansive nor structural stable.

The motivation for our investigation
is the fact that if $U$ is an isolating block of the geometric Lorenz
attractor with vector field $X$ then
for every $Y$ close to $X$ the set
$
\{x\in U:
\omega_Y(x)=\cap_{t\geq 0}Y_t(U)\}
$ is residual in $U$
(this is precisely the same property of the hyperbolic attractors
reported before).
It is then natural to believe
that such a conclusion holds if
$U$ is an isolating block of a singular-hyperbolic attractor.
The answer however is negative as the example \cite[Appendix]{mpu} shows.
Despite we shall
prove that if $U$ is the isolating block of a singular-hyperbolic attractor of $X$, then
the following alternative property holds:
For every vector field $Y$ $C^r$ close to
$X$ the set
$$
\{x\in U:\omega_Y(x) \mbox{ contains a singularity}\}
$$
is {\em residual} in $U$.
In other words, the
positive orbits in a residual subset of $U$ look to be "attracted" to the singularities of $Y$ in $U$.
This fact
can be observed with the computer
in the classical polynomial Lorenz equation \cite{l}.
It
contrasts with the fact that the union of the
stable manifolds of the singularities of $Y$ in $U$
{\em is not residual in any open set}.
We use this property
to prove the persistence
singular-hyperbolic attractors
with only one singularity
as chain-transitive
Lyapunov stable sets.

Now we state our result in a precise way.
Hereafter $M$ denotes a compact Riemannian three manifold unless otherwise stated.
If $U\subset M$
we say that $R\subset U$
{\em residual} if it realizes
as a countable intersection of open-dense subsets of
$U$.
It is well known that
every residual subset
of $U$ is dense in $U$.
Let $X$ be a $C^r$ vector field in $M$ and
let $X_t$ be the
flow generated by $X$, $t\in I\!\! R$.
A compact invariant set
is {\em singular} if it contains a singularity.

\begin{defi}[Attractor]
\label{attractor}
An {\em attracting set} of $X$ is
a compact, invariant, non-empty, set of $X$
equals to $\cap_{t>0}X_t(U)$
for some compact neighborhood $U$ of it.
This neighborhood is called {\em isolating block}.
An {\em attractor} is a
transitive attracting set.
\end{defi}

\begin{rk}
\label{re-1}
\cite{hu} calls attractor what we call
attracting set.
Several definitions
of attractor are considered in \cite{mi}.
\end{rk}

Denote by $m(L)$  and $Det(L)$ the minimum norm
and the Jacobian of a
linear operator $L$ respectively.

\begin{defi}
A compact invariant set $\Lambda$ of
$X$ is
{\it partially hyperbolic}
if there is a continuous invariant
tangent bundle decomposition
$T_\Lambda M=E^s\oplus E^c$ and positive constants
$K,\lambda$ such that

\item{1.}
$E^s$ is contracting:
$
\mid \mid DX_t(x)/E^s_x\mid\mid\leq K e^{-\lambda t}
$,
for every $\forall t>0$ and $x\in \Lambda$;
\item{2.}
$E^s$ dominates $E^c$:
$
\frac{\mid\mid DX_t(x)/E^s_x\mid\mid}{m(DX_t(x)/E^c_x)}
\leq K e^{-\lambda t}
$,
for every $\forall t>0$
and $\forall x\in \Lambda$.

We say that $\Lambda$ has
{\em volume expanding central direction}
if
$$
\mid Det(DX_t(x)/E^c_x)\mid \geq K^{-1}e^{\lambda t},
$$
for every $t>0$ and $x\in \Lambda$.
\end{defi}

A singularity $\sigma$ of $X$ is {\em hyperbolic} if its eigenvalues are not purely
imaginary complex number.

\begin{defi}[Singular-hyperbolic set]
\label{robust}
A compact invariant set of a vector field $X$
is {\em singular-hyperbolic}
if it has singularities
(all hyperbolic) and is partially hyperbolic
with volume
expanding central direction \cite{mpp1}.
A {\em singular-hyperbolic attractor}
is an attractor which is also a singular-hyperbolic set.
\end{defi}

Singular-hyperbolic attractors
cannot be hyperbolic and the most representative example
is the geometric Lorenz \cite{gw}.
Our result is the following.

\begin{mainthm}
\label{th1'}
Let $U$ be an isolating block of a singular-hyperbolic attractor of $X$. If $Y$ is a vector field
$C^r$ close to $X$, then
$\{x\in U:\omega_Y(x)$ is singular$\}$ is residual in $U$
\end{mainthm}

This result is used to prove

\begin{mainthm}
\label{tm'}
Singular-hyperbolic attractors
with only one singularity in $M$
are {\em persistent}
as chain-transitive Lyapunov stable sets.
\end{mainthm}

The precise statement of Theorem \ref{tm'}
(including the definition of chain transitive set, Lyapunov stable set and persistence)
will be given in Section 7.
This paper is organized as follows.
In Section 2 we give some preliminary lemmas.
In particular, Lemma
\ref{attrac'} introduces the {\em continuation}
$A_Y$ of an attracting set $A$ for
nearby vector fields $Y$. In
Definition \ref{stable set} we define
{\em the region of weak attraction} $A_w(Z,C)$ of $C$, where $C$ is a compact invariant sets of a vector field,
as the set of points $z$ such that
$\omega_Z(z)\cap C\neq\emptyset$.
Lemma \ref{dense->residual}
proves that if
$U$ is a neighborhood of $C$ and $A_w(Z,C)\cap U$ is dense
in $U$, then $A_w(Z,C)\cap U$ is residual in $U$.
We finish this section with some
elementary properties of the hyperbolic sets.
We
present two elementary properties of singular-hyperbolic attracting sets in Section 3.

In Section 4 we introduce {\em the Property (P)} for
compact invariant sets $C$ all of whose closed orbits
are hyperbolic.
It requires that
the unstable manifold of every closed orbit
in $C$ intersect transversely the stable manifold
of a singularity in $C$.
This property has been proved for all singular-hyperbolic attractors $\Lambda$
in \cite{mpa1}.
In Lemma \ref{explosion} we prove
that it is open, namely it holds
for the continuation $\Lambda_Y$ of $\Lambda$.
The proof is similar to the one in \cite{mpa1}.

In Section 5 we study the topological dimension
\cite{hw}
of the omega-limit sets
in an isolating block $U$ of a singular-hyperbolic
attracting set with the Property (P).
In particular,
Theorem \ref{th2} proves that
if $x\in U$ then the omega-limit set of $x$ either
contains a singularity or has topological dimension one
provided the stable manifolds of the singularities
in $U$ do not intersect a neighborhood
of $x$.
The proof uses the methods
in \cite{m1} with the Property (P) playing the role
of the transitivity. We need this theorem to apply the
Bowen's theory of one-dimensional hyperbolic sets
\cite{bo}.

In Section 6 we prove Theorem \ref{th1'}.
The proof is based on Theorem \ref{dense} where it is proved that if $U$ is an isolating block
of a singular-hyperbolic attracting set
with the Property (P) of a vector field $Y$,
then $A_w(Y,Sing(Y,U))\cap U$ is dense in $U$ (here $Sing(Y,U)$ denotes the set of singularities of $Y$ in $U$).
The proof follows
applying the Bowen's theory
(that can be used by Theorem \ref{th2})
and the arguments in \cite[p. 371]{mpa1}.
It will follow from Lemma \ref{dense->residual}
applied to $C=Sing(Y,U)$ that
$A_w(Y,Sing(Y,U))\cap U$ is residual in
$U$.
Theorem \ref{th1'}
follows because
$\omega_Y(x)$ is singular
$\forall x\in A_w(Y,Sing(Y,U))\cap U$.
In Section 7 we prove
Theorem \ref{tm'} (see Theorem \ref{tm}).

\section{Preliminary lemmas}

We state some preliminary results.
The first one claims a sort of stability
of the attracting sets. It seems to be well
known and we prove it here for completeness.
If $M$ is a manifold and $U\subset M$ we denote by
$int(U)$ and $clos(U)$ the interior
and the closure of $U$ respectively.

\begin{lemma}[Continuation of attracting sets]
\label{attrac'}
Let $A$ be an attracting set
{\em containing a hyperbolic closed orbit}
of a $C^r$ vector field $X$.
If $U$ is an isolating block of $A$, then
for every vector field
$Y$ $C^r$ close to $X$ {\em the continuation}
$$
A_Y=\cap_{t\geq 0}Y_t(U)
$$
{\em of $A$
in $U$}
is an attracting set with isolating block $U$ of $Y$.
\end{lemma}

\begin{proof}
Since $A$ contains a hyperbolic closed orbit
we have that $A_Y\neq\emptyset$ for every $Y$ close to $X$
(use for instance the Hartman-Grobman
Theorem \cite{dmp}). Since $U$ is compact we have that
$A_Y$ also does.
Then, to prove the lemma,
we only need to prove that
if $Y$ is close to $X$ then $U$ is a
compact neighborhood of $A_Y$.
For this we proceed as follows.
Fix an open set $D$ such that
$$
A\subset D\subset clos(D)\subset int(U)
$$
and for all $n\in I\!\! N$ we
define
$$
U_n=\cap_{t\in [0,n]}X_t(U).
$$
Clearly $U_n$ is a compact set sequence
which is nested ($U_{n+1}\subset U_n$) and satisfies
$
A=\cap_{n\in I\!\! N}U_n
$.
Because $U_n$ is nested we
can find $n_0$ such that
$U_{n_0}\subset D$.
In other words
$$
\cap_{t\in [0,n_0]}X_t(U)\subset D.
$$
Taking complement one has
$$
M\setminus D
\subset \cup_{t\in [0,n_0]}X_t(M\setminus U).
$$
But $X_t(M\setminus U)$ is open
($\forall t$) since $U$ is compact and $X_t$ is a diffeomorphism. Hence
$\{X_t(M\setminus U):t\in [0,n_0]\}$ is an open covering
of $M\setminus D$. Because $D$ is open we
have that
$M\setminus D$ is compact and so
there are finitely many
$t_1,\cdots ,t_k\in [0,n_0]$ such that
$$
M\setminus D\subset
X_{t_1}(M\setminus U)\cup\cdots
\cup X_{t_k}(M\setminus U).
$$
By the continuous dependence of $Y_t(U)$
on $Y$ (with $t$ fixed)
one has
$$
M\setminus D\subset
Y_{t_1}(M\setminus U)\cup\cdots
\cup Y_{t_k}(M\setminus U)
$$
for all $Y$ $C^r$ close to $X$.
By taking complement once more we obtain
$$
Y_{t_1}(U)\cap \cdots \cap Y_{t_k}(U)\subset D.
$$
As $t_1,\cdots ,t_k\geq 0$ one has
$\cap_{t\in [0,n_0]}Y_t(U)\subset
Y_{t_1}(U)\cap \cdots \cap Y_{t_k}(U)$ and then
$$
\cap_{t\in [0,n_0]}Y_t(U)\subset D
$$
for every $Y$ close to $X$.
On the other hand,
it follows from the definition that
$A_Y\subset \cap_{t\in [0,n_0]}Y_t(U)$ and so
$A_Y\subset D$
for every $Y$ close to $X$.
Because $clos(D)\subset int(U)$
we have that $A_Y\subset int(U)$.
This proves that $U$ is a compact neighborhood
of $A_Y$ and the lemma follows.
\end{proof}

\begin{rk}
\label{closeness}
The above proof shows
that the compact set-valued map
$Y\to A_Y$ is continuous in the following sense:
For every open set $D$ containing
$A$ one has $A_Y\subset D$ for every $Y$ $C^r$ close to $X$.
Such a continuity
is weaker than the continuity
with respect to the
Hausdorff metric.
It follows from the above-mentioned continuity that
if $A$ is a singular-hyperbolic attracting set
of $X$ and $Y$ is close to $X$,
then the continuation
$A_Y$ in $U$ is
a singular-hyperbolic attracting set of
$Y$.
\end{rk}

The following definition can be found in
\cite[Chapter V]{bs}.

\begin{defi}[Region of attraction]
\label{stable set}
Let $C$ be a compact invariant set of
a vector field $Z$.
We define
{\em the region of attraction}
and {\em the region of weak attraction} of $C$
by 
$$
A(C)=\{x\in M:
\omega_X(p)\subset C\}
\,\,\,\,\mbox{and}\,\,\,\,
A_w(C)=
\{z:\omega_Y(z)\cap C\neq\emptyset\}
$$
respectively.
We shall write $A(Z,C)$ and $A_w(Z,C)$ to indicate dependence on $Z$.
\end{defi}

The region of attraction is also
called {\em stable set}.
The inclusion below is obvious
\begin{equation}
\label{paja}
A(Z,C)
\subset A_w(Z,C).
\end{equation}

The elementary lemma below will be used
in Section 6. Again we prove it for the sake of completeness.

\begin{lemma}
\label{dense->residual}
If $C$ a compact invariant set
of a vector field $Z$ and $U$ is
a compact neighborhood of $C$, then
the following properties are equivalent:
\item{1.}
$A_w(Z,C)\cap U$ is dense in $U$
\item{2.}
$A_w(Z,C)\cap U$ is residual in $U$.
\end{lemma}

\begin{proof}
Clearly (2) implies (1).
Now we assume (1) namely
$A_w(Z,C)\cap U$ is dense in $U$.
Defining
$$
W_n=
\{x\in U:Z_t(x)\in B_{1/n}(C)
\mbox{ for some }t>n\} \,\,\,\,\forall n\in I\!\! N
$$
one has
$$
A_w(Z,C)\cap U=\cap_n W_n.
$$
In particular $A_w(Z,C)\cap U\subset W_n$ for all $n$.
Hence $W_n$ is dense in $U$
(for all $n$) since
$A_w(Z,C)\cap U$ does.
On the one hand, $W_n$ is open in $U$
\cite[Tubular
Flow-Box Theorem]{dmp} because
$B_{1/n}(T)$ is open.
This proves that $W_n$ is open-dense in $U$
and the result follows.
\end{proof}

Next we state the classical definition of hyperbolic set.

\begin{defi}[Hyperbolic set]
\label{hyperbolic-set}
A compact, invariant set $H$ of
a $C^1$ vector field $X$ is
{\em hyperbolic} if
there are a continuous, tangent bundle,
invariant,
splitting $T\Lambda=E^s\oplus E^X\oplus E^u$
and positive constants $C,\lambda$
such that
$\forall x\in H$ one has:

\item{1.}
$E^X_x$ is the direction of $X(x)$ in $T_xM$.
\item{2.}
{\it $E^s$ is contracting}:
$
\mid\mid DX_t(x)/E^s_x\mid\mid\leq
Ce^{-\lambda t}$, $\forall t\geq 0$.
\item{3.}
{\it $E^u$ is expanding}:
$
\mid\mid DX_{t}(x)/E^u_{x}\mid\mid
\geq
C^{-1}e^{\lambda t}$, $\forall t\geq 0$.

A closed orbit
of $X$ is hyperbolic if
it is hyperbolic as a compact, invariant set of $X$.
A hyperbolic set is {\it saddle-type}
if $E^s\neq 0$ and $E^u\neq 0$.
\end{defi}

The Invariant Manifold Theory
\cite{hps} says that
through each point $x\in H$
pass smooth injectively immersed submanifolds
$W^{ss}(x),W^{uu}(x)$ tangent to $E^s_x,E^u_x$ at
$x$.
The manifold $W^{ss}(x)$,
the strong stable manifold at $x$, is characterized by
$y\in W^{ss}(x)$ if and only
if $d(X_t(y),X_t(y))$ goes to $0$ exponentially as
$t\to\infty$.
Similarly $W^{uu}(x)$,
the strong unstable manifold at $x$,
is characterized by $y\in W^{uu}(x)$
if and only if $d(X_t(y),X_t(x))$ goes to $0$ exponentially as $t\to -\infty$.
These manifolds
are invariant, i.e.
$
X_t(W^{ss}(x))=
W^{ss}(X_t(x))$ and
$X_{t}(W^{uu}(x))=
W^{uu}(X_{t}(x))$,
$\forall t$.
For all $x,x'\in H$ we have that
$W^{ss}(x)$ and $W^{ss}(x')$
either coincides or are disjoint.
The maps
$x\in H\to W^{ss}(x)$
and $x\in H\to W^{uu}(x)$
are continuous (in compact parts).
For all $x\in H$ we define
$$
W^s_X(x)=\cup_{t\in I\!\! R}W^{ss}(X_t(x))\,\,\,\,
\mbox{and}\,\,\,\,
W^u_X(x)=
\cup_{t\in I\!\! R}W^{uu}(X_t(x)).
$$
Note that if $O\subset H$ is a closed orbit
then
$$
A(X,O)=W^s_X(O)
$$
but $A_w(X,O)\neq W^s_X(O)$ in general.
If $H$ is saddle-type and $dim(M)=3$, then
both $W^s_X(x),W^u_X(x)$ are one-dimensional submanifolds
of $M$. In this case given $\epsilon>0$ we denote by
$W^{ss}_X(x,\epsilon)$ an interval of length
$\epsilon$ in $W^{ss}_X(x)$ centered at $x$
(this interval is often called the local
strong stable manifold of $x$).

\begin{defi}
Let $\{O_n:n\in I\!\! N\}$ be a sequence of hyperbolic periodic
orbits of $X$. We say that
{\em the size of $W^s_X(O_n)$ is uniformly bounded
away from zero} if
there is $\epsilon>0$
such that the
local strong stable manifold
$W^{ss}_X(x_n,\epsilon)$ is well defined
for every $x_n\in O_n$ and every $n\in I\!\! N$.
\end{defi}

\begin{rk}
\label{strong-stable}
Let $O_n$ be a sequence of hyperbolic periodic orbits of
a vector field $X$.
It follows from the Stable Manifold Theorem for hyperbolic sets
\cite{hps} that
the size of $W^s_X(O_n)$ is uniformly bounded
away from zero if all the periodic
orbits $O_n$ ($n\in I\!\! N$) are contained
in the same hyperbolic set $H$ of $X$.
\end{rk}

\section{Two lemmas for singular-hyperbolic attracting sets}

Hereafter we denote by
$M$ a compact three manifold.
Recall that $clos(\cdot)$ denotes the closure
of $(\cdot)$.
In addition,
$B_\delta(x)$ denotes
the (open)
$\delta$-ball in $M$ centered at $x$.
If $H\subset M$ we
denote
$B_\delta(H)=\cup_{x\in H}B_\delta (x)$. 
For every vector field $X$ on $M$ we
denote by $Sing(X)$ the set of singularities
of $X$ and if $B\subset M$ we define
$Sing(X,B)=Sing(X)\cap B$.

\begin{lemma}
\label{hyperbolic}
Let $\Lambda$ be a singular-hyperbolic
attracting set of a $C^r$ vector field $Z$ on $M$.
Let $U$ be an isolating block of $\Lambda$.
If $x\in U$ and $\omega_Z(x)$ is non-singular,
then every $k\in \omega_Z(x)$ is
accumulated by a hyperbolic periodic orbit sequence $\{O_n:n\in I\!\! N\}$ such that the size of $W^s_Z(O_n)$ is
uniformly bounded away from zero.
\end{lemma}

\begin{proof}
For every $\epsilon>0$ we define
$$
\Lambda_\epsilon=
\cap_{t\in I\!\! R}Z_t(\Lambda\setminus
B_\epsilon(Sing(Z,\Lambda)).
$$
Clearly $\Lambda_\epsilon$ is either $\emptyset$
or a compact, invariant, non-singular set of
$Z$. If $\Lambda_\epsilon\neq\emptyset$,
then $\Lambda_\epsilon$ is hyperbolic \cite{mpp2}.
Observe that $\omega_X(x)$ is non-singular by assumption.
Then, there are $\epsilon>0$ and $T>0$ such that
$$
Z_t(x)\notin clos(B_\epsilon(Sing(Z,U))),
\,\,\,\,\,\,\,\,\forall t\geq T.
$$
It follows that
$\omega_Z(x)\subset \Lambda_\epsilon$
and so $\Lambda_\epsilon\neq\emptyset$ is a hyperbolic set.
In addition,
for every $\delta>0$ there is
$T_\delta>0$ such that
$$
Z_t(x)\in B_\delta(\Lambda_\epsilon),
$$
for every $t>T_\delta$.
Pick $k\in \omega_Z(x)$. The last
property implies that
for every $\delta>0$ there is
a periodic $\delta$-pseudo-orbit
in $B_\delta(\Lambda_\epsilon))$
formed by paths in the positive $Z$-orbit
of $x$.
Applying the Shadowing Lemma for Flows
\cite[Theorem 18.1.6 pp. 569]{hk}
to the hyperbolic set $\Lambda_\epsilon$
we arrange a periodic orbit sequence
$O_n\subset \Lambda_{\epsilon/2}$ accumulating
$k$.
Then, Remark \ref{strong-stable} applies
since $H=\Lambda_{\epsilon/2}$ is hyperbolic and
contains $O_n$ (for all $n$).
The lemma is proved.
\end{proof}

The following is a minor modification
of \cite[Theorem A]{m2}.

\begin{lemma}
\label{explo}
If $U$ is an isolating block of a singular-hyperbolic attractor
of a $C^r$ vector field $X$ in $M$,
then every attractor in $U$ of
every vector field $C^r$ close to $X$ is
singular.
\end{lemma}

\begin{proof}
Let $\Lambda$ be the singular-hyperbolic attractor of $X$ having $U$ as isolating block.
By \cite[Theorem A]{m2}
there is a neighborhood
$D$ of $\Lambda$ such that
every attractor of every
vector field $Y$ $C^r$ close to $X$ is singular.
By Remark
\ref{closeness} we have that $
\cap_{t\geq 0}Y_t(U)\subset D$ for all $Y$ close to $X$.
Now if $A\subset U$ is an attractor of $Y$, then $A\subset \cap_{t\geq 0}Y_t(U)$
by invariance.
We conclude that $A\subset D$ and then $A$ is singular
for all $Y$ close to $X$.
This proves the lemma.
\end{proof}

\section{Property (P)}

First we state the definition.
As usual we write
$S\pitchfork S'\neq\emptyset$
to indicate that there is a transverse intersection point between
the submanifolds $S,S'$.

\begin{defi}[The Property (P)]
\label{main-lemma}
Let $\Lambda$ be a compact invariant set of
a vector field $X$. Suppose that
all the closed orbits of $\Lambda$ are hyperbolic.
We say that $\Lambda$ satisfies {\em the Property (P)}
if for every point $p$ on a periodic orbit of $\Lambda$
there is $\sigma\in Sing(X,\Lambda)$ such that
$$
W^u_Y(p)\pitchfork W^s_Y(\sigma)\neq\emptyset.
$$
\end{defi}

The lemma below is a direct consequence
of the classical Inclination-lemma \cite{dmp}
and the transverse intersection in Property (P).

\begin{lemma}
\label{inclination-lemma}
Let $\Lambda$ a compact invariant set with
the Property (P) of a vector field $Z$ in a
manifold $M$ and $I$ be a submanifold of $M$.
If
there is a periodic orbit
$O\subset \Lambda$ of $Z$ such that
$$
I\pitchfork W^s_Z(O)\neq\emptyset,
$$
then
$$
I\cap\left(
\cup_{\sigma\in Sing(Z,\Lambda)}W^s_Z(\sigma)\right)\neq\emptyset.
$$
\end{lemma}

The Property (P) was proved
in \cite[Theorem 5.1]{mpa1} for
all singular-hyperbolic attractors.
Here we prove that such a property is open, namely
it holds for the continuation in Lemma \ref{attrac'} of a singular-hyperbolic attractor.

\begin{lemma}[Openness of the Property (P)]
\label{explosion}
Let $U$ be an isolating block of a singular-hyperbolic attractor
of a $C^r$ vector field $X$ on $M$. Then,
the continuation
$$
\Lambda_Y=\cap_{t\geq 0}Y_t(U)
$$
has the Property (P)
for every vector field $Y$ $C^r$ close to $X$.
\end{lemma}

\begin{proof}
By Lemma \ref{attrac'}
we have that $\Lambda_Y$ is an attracting set with isolating block $U$ since
$\Lambda$ has a hyperbolic singularity.
Now let $p$ be a point of a periodic
orbit $\gamma\subset \Lambda_Y$ of $Y$.
Then
$$
clos(W^u_Y(p))\subset \Lambda_Y
$$
since $\Lambda_Y$ is attracting.
We claim
$$
clos(W^u_Y(p))\cap Sing(Y,U)\neq\emptyset.
$$
Indeed suppose that it is not so, i.e.
there is $Y$ $C^r$ close to $X$ such
that $clos(W^u_Y(p))\cap Sing(Y,U)=\emptyset$
for some $p$ in a periodic orbit of $Y$ in $U$.
It follows from \cite{mpp2} that
$clos(W^u_Y(p))$ is a hyperbolic set.
Since $W^u_Y(p)$ is a two-dimensional
submanifold we can easily prove that $clos(W^u_Y(p))$
is an attracting set of $Y$.
This attracting set necessarily contains
a hyperbolic attractor $A$ of $Y$.
Since $A\subset clos(W^u_Y(p))\subset \Lambda_Y
\subset U$ we conclude that
$A\subset U$.
By Lemma \ref{explo} we have that
$A$ is singular as well.
We conclude that $A$ is an attracting singularity
of $Y$ in $U$. This contradicts
the volume expanding condition
at Definition \ref{robust} and the claim follows.
One completes the proof
of the lemma using the claim
as in \cite[Theorem 5.1]{mpa1}.
\end{proof}

\section{Topological dimension and the Property (P)}

We study the topological dimension
of the omega-limit set in an isolating block
of a singular-hyperbolic attracting set with the Property (P).
First of all we recall the
classical definition of topological dimension
\cite{hw}.

\begin{defi}
\label{top-dim}
The {\em topological dimension} of a
space $E$ is either $-1$
(if $E=\emptyset$) or
the last integer
$k$ for which every point has arbitrarily small
neighborhoods whose boundaries
have dimension less than $k$.
A space with topological dimension
$k$ is said to be {\em k-dimensional}.
\end{defi}

The result of this section is the following.

\begin{thm}
\label{th2}
Let $U$ be an isolating block of a singular-hyperbolic attracting set with the Property (P) of a $C^r$ vector field $Y$ on $M$.
If $x\in U$ and there is
$\delta>0$ such that
$$
B_\delta(x)\cap \left(
\cup_{\sigma\in Sing(Y,U)}W^s_Y(\sigma)\right)=
\emptyset ,
$$
then $\omega_Y(x)$ is either singular or a one-dimensional hyperbolic set.
\end{thm}

\begin{proof}
Let $\Lambda_Y$ be the singular-hyperbolic attracting set
of $Y$ having $U$ as isolating block.
Obviously $Sing(Y,U)=Sing(Y,\Lambda_Y)$.
Let $x,\delta$ be
as in the statement. Define
$$
H=\omega_Y(x).
$$
We shall assume that $H$ is non-singular.
Then $H$ is a hyperbolic set
by \cite{mpp2}.
To prove that $H$
is one-dimensional we shall
use the arguments in \cite{m1}.
However we have to take some
care because
$\Lambda$ is not transitive.
The Property (P) will supply an alternative argument.
Let us present the details.

First we note that
by Lemma \ref{hyperbolic} every point $k\in H$ is accumulated
by a periodic orbit sequence
$O_n$ satisfying the conclusion of that lemma.
Second, by the Invariant Manifold Theory
\cite{hps},
there is an invariant contracting foliation
$\{{\mathcal F}^s(w):w\in \Lambda_Y\}$ which is tangent
to the contracting direction of $Y$ in $\Lambda_Y$.
A cross-section of $Y$ will be a $2$-disk transverse to $Y$.
When $w\in \Lambda_Y$ belongs to a $2$-disk
$D$ transverse to $Y$,
we define
${\mathcal F}^s(w,D)$ as the connected component
containing $w$ of
the projection of
${\mathcal F}^s(w)$ onto $D$ along the flow of $Y$.
The boundary and the interior
of $D$ (as a submanifold of $M$) are
denoted by $\partial D$ and $int(D)$ respectively.
$D$ is
a {\em rectangle} if it is diffeomorphic to
the square $[0,1]\times [0,1]$.
In this case $\partial D$
as a submanifold of $M$ is formed
by four curves $D^t_h, D^b_h, D^l_v,D^r_v$
($v$ for vertical, $h$ for horizontal,
$l$ for left, $r$ for right, $t$ for top
and $b$ for bottom).
One defines vertical and horizontal
curves in $D$ in the natural way.

Now we prove a sequence of lemmas corresponding
to lemmas 1-4 in \cite{m1} respectively.

\begin{lemma}
\label{claim0}
For every regular point $z\in \Lambda_Y$ of $Y$
there is a rectangle $\Sigma$ such that the properties below hold:
\item{1.}
$z\in int(\Sigma)$;
\item{2.}
If $w\in \Lambda_Y$ then ${\mathcal F}^s(w,\Sigma)$ is
a horizontal curve in $\Sigma$;
\item{3.}
If $\Lambda_Y\cap \Sigma^t_h\neq\emptyset$ then
$\Sigma^t_h={\mathcal F}^s(w,\Sigma)$ for some
$w\in \Lambda_Y\cap \Sigma$;
\item{4.}
If $\Lambda_Y\cap \Sigma^b_h$ then
$\Sigma^b_h={\mathcal F}^s(w,\Sigma)$ for some
$w\in \Lambda_Y\cap \Sigma$.
\end{lemma}

\begin{proof}
The proof of this lemma is
similar to \cite[Lemma 1]{m1}.
Observe that the corresponding
proof in \cite{m1} does not use
the transitivity hypothesis.
\end{proof}

\begin{defi}
If $w\in H\cap \Sigma$ we
denote by
$(H\cap \Sigma)_w$ the connected component
of $H\cap \Sigma$ containing
$w$.
\end{defi}

With this definition we shall
prove the following lemma.

\begin{lemma}
\label{claim1}
If
$w\in H\cap \Sigma$ and
$(H\cap \Sigma)_w\neq\{w\}$, then
$(H\cap \Sigma)_w$ contains
a non-trivial curve in
the union
${\mathcal F}^s(w,\Sigma)\cup \partial \Sigma$.
\end{lemma}

\begin{proof}
We follow the same steps of the proof
of Lemma 2 in \cite{m1}.
First we observe that
$
(H\cap \Sigma)_x\cap
(int(\Sigma)\setminus
{\mathcal F}^s(x,\Sigma))\neq\emptyset
$.
Hence we can fix $w'\in (H\cap \Sigma)_x\cap
(int(\Sigma)\setminus {\mathcal F}^s(x,\Sigma))$.
Clearly ${\mathcal F}^s(w',\Sigma)$ is a horizontal curve
which together with ${\mathcal F}^s(w,\Sigma)$ form
the horizontal boundary curves of a rectangle $R$ in $\Sigma$.
One has that $H\cap int(B)\neq\emptyset$
for, otherwise, $w$ and $w'$ would be in different
connected components of $H\cap \Sigma$ a contradiction.
Hence we can choose
$h\in H\cap int(B)$.
Since $H=\omega_Y(y)$ we have that
there is $y'$ in the positive $Y$-orbit of $y$
arbitrarily close to $h$.
In particular, $y'\in int(B)$.
By the continuity of the foliation ${\mathcal F}^s$
we have that ${\mathcal F}^s(y',\Sigma)$
is a horizontal curve separating $\Sigma$ in two
connected components containing $w$ and $w'$ respectively.
Since $w,w'$ belong to the same connected component
of $H\cap \Sigma$ we conclude that there is
$k\in {\mathcal F}^s(y',\Sigma)\cap H\neq\emptyset$.

On one hand, by Lemma \ref{hyperbolic},
$k\in H$ is accumulated
by a hyperbolic periodic orbit sequence $O_n$ such that
the size of $W^s_Y(O_n)$ is uniform
bounded away from zero.
On the other hand
$y'$ belongs to the positive orbit
of $y$ and $y\in B_\delta(x)$.
By the uniform
size of $W^s_Y(O_n)$
one has
$
B_\delta(x)\cap W^s_Y(O_n)\neq\emptyset
$
for some $n\in I\!\! N$.
Since $B_\delta(x)$ is open we conclude that
$$
B_\delta(x)\pitchfork W^s_Y(O_n)\neq\emptyset
$$
Then,
$$
B_\delta(x)\cap\left(\cup_{\sigma\in Sing(Y,U)}W^s_Y(\sigma)\right)\neq\emptyset
$$
by
Lemma \ref{inclination-lemma}
since $\Lambda_Y$ has the Property (P).
This is a contradiction which proves the lemma.
\end{proof}

\begin{lemma}
\label{claim2}
For every $w\in H$ there is a rectangle
$\Sigma_w$ containing $w$ in its interior
such that
$H\cap \Sigma_w$ is $0$-dimensional.
\end{lemma}

\begin{proof}
This lemma corresponds to Lemma 3 in \cite{m1}
with similar proof.
Let $\Sigma_w=\Sigma$ where
$\Sigma$ is given by Lemma \ref{claim1}.
Let $J\subset
{\mathcal F}^s(w,\Sigma)\cap \partial \Sigma$ be the curve in the conclusion of this lemma.
We can assume that
$J$ is contained in either
${\mathcal F}^s(w,\Sigma)$ or $\partial \Sigma$.
If $J\subset {\mathcal F}^s(w,\Sigma)$
we can prove as in the proof of [M3, Lemma 3]
that $y\in H$ and so
$y$ is accumulated by periodic
orbits whose unstable and stable manifolds
have uniform size.
We arrive a contradiction by Lemma \ref{explosion} as
in the last part of the proof of Lemma \ref{claim1}.
Hence we can assume that
$J\subset \partial \Sigma$.
We can further assume that $J\subset \Sigma_v^l$
(say) for otherwise we get a contradiction
as in the previous case.
Now if
$J\subset \Sigma_v^l$
then we can obtain a contradiction as before
again using the Property (P) and
Lemma \ref{inclination-lemma}.
This proves the result.
\end{proof}

The following lemma corresponds
to \cite[Lemma 4]{m1}.

\begin{lemma}
\label{claim3}
$H$ can be covered by a finite collection of
closed one-dimensional subsets.
\end{lemma}

\begin{proof}
If $w\in H$ we consider
the cross-section $\Sigma_w$
in Lemma \ref{claim3}. By saturating
forward and backward $\Sigma_w$ by the flow of $Y$ we
obtain a compact neighborhood of $w$
which is one-dimensional
(see \cite[Theorem III 4 p. 33]{hw}).
Hence there is a neighborhood covering
of $H$ by compact one-dimensional sets.
Such a covering has a finite subcovering
since $H$ is compact. Such a subcovering
proves the result.
\end{proof}

Theorem \ref{th2} now follows from Lemma \ref{claim3}
and \cite[Theorem III 2 p. 30]{hw}.
\end{proof}

\section{Proof of Theorem \ref{th1'}}

The proof is based on the following result.

\begin{thm}
\label{dense}
Let $U$ be an isolating block of a singular-hyperbolic attracting set with the Property (P) of a vector field $Y$ on $M$. Then
$A_w(Y,Sing(Y,U))\cap U$ is residual in $U$.
\end{thm}

\begin{proof}
By Lemma \ref{dense->residual} it suffices
to prove that $A_w(Y,Sing(Y,U))\cap U$ is dense in $U$.
Let $\Lambda_Y$ be the singular-hyperbolic attracting set
of $Y$ having $U$ as isolating block.
Obviously $Sing(Y,U)=Sing(Y,\Lambda_Y)$.
To simplify the notation we
write $R_Y=A_w(Y,Sing(Y,U))\cap U$.
Suppose by contradiction
that $R_Y$ is not dense in $U$.
Then, there is $x\in U$ and $\delta>0$ such that
$B_\delta(x)\cap R_Y=\emptyset$.
In particular,
$\omega_Y(x)\cap Sing(Y,U)=\emptyset$ and so
$\omega_Y(x)$ is non-singular.
Recalling the inclusion Eq.(\ref{paja}) at Section 2 one has
$$
U\cap\left(\cup_{\sigma\in Sing(Y,U)}W^s_Y(\sigma)\right)\subset
R_Y.
$$
Thus
\begin{equation}
\label{absurd}
B_\delta(x)\cap \left(
\cup_{\sigma\in Sing(Y,U)}W^s_Y(\sigma)\right)=
\emptyset.
\end{equation}
It then follows from Theorem
\ref{th2} that
$H=\omega_Y(x)$ is a one-dimensional hyperbolic set.
This allows to apply the Bowen's
Theory \cite{bo} of one-dimensional hyperbolic sets.
More precisely
there is
a family of (disjoint) cross-sections
${\mathcal S}=\{S_1,\cdots ,S_r\}$ of small diameter
such that
$H$ is the flow-saturated
of $H\cap int({\mathcal S}')$, where
$
{\mathcal S}'=\cup S_i$
and $int({\mathcal S}')$ denotes the interior
of ${\mathcal S}'$ (as a submanifold).
Next we choose an interval $I$
tangent to the central direction $E^c$ of $Y$
in $U$ such that
$$
x\in I\subset B_\delta(x).
$$
We choose $I$ to be transverse to the direction
$E^Y$ induced by $Y$.
Since
$E^c$ is volume expanding and
$H$ is non-singular we have that
the Poincar\'e map
induced by $X$ on ${\mathcal S}'$
is expanding along $I$.
As
in \cite[p. 371]{mpa1}
we can find $\delta'>0$ and a open arc sequence
$J_n\subset {\mathcal S}'$ in the positive orbit
of $I$ with length $\geq \delta'$
such that there is
$x_n$ {\em in the positive orbit of $x$}
contained in the interior of
$J_n$.
We can fix $S=S_i\in {\mathcal S}$
in order to assume that
$J_n\subset S$ for
every $n$.
Let $w\in S$ be a limit point of
$x_n$.
Then $w\in H\cap int({\mathcal S}')$.
Because $I$ is tangent to
$E^c$ the interval sequence $J_n$
converges to an interval
$J\subset W^u_Y(w)$ in the $C^1$ topology
($W^u_Y(w)$ exists because $w\in H$ and $H$
is hyperbolic).
$J$ is not trivial since the length
of $J_n$ is $\geq \delta'$.
It follows from this lower bound that
$J_n$ intersects $W^s_Y(w)$ for some
$n$ large.
Now, by Lemma \ref{hyperbolic}, $w$ is accumulated by
periodic orbits $O_n$
satisfying the conclusion of this lemma.
The continuous dependence in compact parts
of the stable manifolds implies
$
J_n\pitchfork W^s_Y(O_n)\neq\emptyset
$.
Since $J_n$ is in the positive orbit of $I$ and
$I\subset B_\delta(x)$ we obtain
$$
B_\delta(x)\pitchfork W^s_Y(O_n)\neq\emptyset.
$$
Then,
$$
B_\delta(x)\cap \left(\cap_{\sigma\in Sing(Y,U)}
W^s_Y(\sigma)\right)\neq\emptyset
$$
by Lemma \ref{inclination-lemma}
since $\Lambda_Y$ has the Property (P).
This is a contradiction by  Eq.(\ref{absurd}).
This contradiction proves
that $R_Y$ is dense in $U$ for all
$Y$ $C^r$ close to $X$.
\end{proof}

\noindent
{\bf Proof of Theorem \ref{th1'}:}
Let $U$ be an isolating block of a singular-hyperbolic attractor
of a $C^r$ vector field $X$ on $M$.
By Lemma \ref{attrac'} we have that $\Lambda_Y=\cap_{t\geq 0}Y_t(U)$
is a singular-hyperbolic attracting set
with isolating block $U$ for all vector field $Y$ $C^r$ close to $X$.
In addition,
$\Lambda_Y$ has the Property (P)
by Lemma \ref{explosion}.
It follows from Theorem \ref{dense} that
$A_w(Y,Sing(Y,U))\cap U$ is residual in $U$.
The result follows
because
$\omega_Y(x)$ is singular
$\forall x\in A_w(Y,Sing(Y,U))\cap U$
(recall Definition \ref{stable set}).
\qed

\begin{rk}
Let $Y$ be a vector field in a manifold $M$. In
\cite[Chapter V]{bs} it was
defined a {\em weak attractor} of $Y$ as a closed set
$C\subset M$ such that $A_w(Y,C)$ is a neighborhood
of $C$.
Similarly one can define
a {\em generic weak attractor} of $Y$ as a
closed set $C\subset M$
such that $A(Y,C)\cap U$ is residual in $U$
for some neighborhood $U$ of $C$
(compare with the definition of generic attractor
\cite[Appendix 1 p.186]{mi}).
A direct consequence of Theorem \ref{dense} is
that the set of singularities of
a singular-hyperbolic attractor of $Y$
is a generic weak attractor of $Y$.
\end{rk}

\section{Persistence of singular-hyperbolic attractors}

In this section we prove Theorem \ref{tm'}
as an application of Theorem \ref{th1'}.
The idea is to address the question below
which is a weaker local version of
the Palis's conjecture \cite{p}.

\begin{q}
\label{q1}
Let $\Lambda$ an attractor
of a $C^r$ vector field $X$ on $M$
and $U$ be an isolating block of $\Lambda$.
Does every vector field $C^r$ close to $X$
exhibit an attractor in $U$?
\end{q}

This question
has positive answer
for hyperbolic attractors,
the geometric Lorenz attractors and the example in \cite{mpu}.
In general we give a partial positive
answer for all singular-hyperbolic attractors
with only one singularity
in terms of chain-transitive
Lyapunov stable sets.

\begin{defi}
\label{LS}
A compact invariant set $\Lambda$ of
a vector field $X$ is
{\em Lyapunov stable}
if for every
open set $U\supset \Lambda$ there is an open
set $\Lambda\subset V\subset U$ such that
$\cup_{t>0}X_t(V)\subset U$.
\end{defi}

Recall that $B_\delta(x)$\
denotes the (open) ball centered at $x$ with radius $\delta>0$. 

\begin{defi}
\label{chain-tran}
Given $\delta>0$ we define
a {\em $\delta$-chain} of $X$ as a pair
of finite sequences $q_1,...,q_{n+1}\in M$ and
$t_1,...,t_n\geq 1$
such that
$$
X_{t_i}(B_\delta(q_i))\cap B_\delta (q_{i+1})\neq \emptyset, \,\,\,\,\forall i=1,\cdots , n.
$$
The $\delta$-chain joints
$p,q$ if $q_1=q$ and $q_{n+1}=p$.
A compact invariant set
$\Lambda$ of $X$ is
{\em chain-transitive} if
every pair of points $p,q\in \Lambda$
can be joined by a $\delta$-chain, $\forall \delta>0$.
\end{defi}

Every attractor is a chain-transitive Lyapunov stable set but not vice versa.
The following generalizes
the concept of robust transitive attractor
(see for instance \cite{mpa2}).

\begin{defi}
\label{persist}
Let $\Lambda$ be a chain-transitive Lyapunov stable set
of a $C^r$ vector field $X$, $r\geq 1$.
We say that $\Lambda$ is
{\em $C^r$ persistent} if for every
neighborhood $U$ of $\Lambda$ and
every vector field $Y$ $C^r$ close to
$X$ there is a chain-transitive Lyapunov stable
set $\Lambda_Y$ of $Y$ in $U$ such that
$A(Y,\Lambda_Y) \cap U$
is residual in $U$.
\end{defi}

Compare this definition
with the one in \cite{hu} where it is
required the continuity of $Y\to \Lambda_Y$
(with respect to the
Hausdorff metric) instead of the residual condition of the stable set.
Another related definition is
that of $C^r$ weakly robust attracting sets
in \cite{cmp}.
The result of this section is the following one.
It is precisely the Theorem \ref{tm'}
stated in the
Introduction.

\begin{thm}
\label{tm}
Singular-hyperbolic attractors
with only one singularity
for $C^r$ vector fields on $M$
are $C^r$ persistent.
\end{thm}

\begin{proof}
Let $\Lambda$ be a singular-hyperbolic attractor
of a $C^r$ vector field $X$ on $M$.
Suppose that $\Lambda$ contains a unique singularity
$\sigma$.
Let $U$ be a neighborhood
of $\Lambda$. We can suppose that
$U$ is an isolating block.
Let $\sigma(Y)$ the continuation of $\sigma$ for every
vector field $Y$ close to $X$. Note that $\sigma(X)=\sigma$. Clearly $Sing(Y,U)=\{\sigma(Y)\}$
for every $Y$ close to $X$.

For every vector field $Y$ $C^r$ close to
$X$ one defines
$$
\Lambda(Y)=\{
q\in \Lambda_Y:\forall \delta>0\,\,
\exists\mbox{$\delta$-chain joining $\sigma(Y)$ and }
q
\}.
$$
Recall that $\Lambda_Y$ is the continuation
of $\Lambda$ in $U$ for $Y$ close to $X$ as in
Lemma \ref{attrac'}.
We note that $\Lambda(Y)\neq\Lambda_Y$
in general \cite{mpu}.

To prove the theorem we shall
prove that $\Lambda(Y)$ satisfies
the following properties ($\forall Y$ $C^r$ close to $X$):

\item{(1)}
$\Lambda(Y)$ is Lyapunov stable.
\item{(2)}
$\Lambda(Y)$ is chain-transitive.
\item{(3)}
$A(Y,\Lambda(Y))\cap U$
is residual in $U$.

One can easily prove (1).
To prove (2) we
pick $p,q\in \Lambda(Y)$ for $Y$ close to $X$ and
fix $\delta>0$.
By Theorem \ref{th1'} there is
$x\in B_\delta(p)$ such that $\omega_Y(x)$
contains $\sigma(Y)$. Hence
there is $t>1$ such that
$X_t(x)\in B_\delta(\sigma)$.
On the other hand, since $q\in \Lambda(Y)$,
there is a $\delta$-chain
$(\{t_1,\cdots ,t_n\},\{q_1,\cdots ,q_{n+1}\})$
joining $\sigma$ to $q$.
Then (2) follows since
the $\delta$-chain
$(\{t,t_1,\cdots ,t_n\},\{x,q_1,\cdots ,q_{n+1}\})$
joints $p$ and $q$.
To finish we prove (3).
It follows from well known properties
of Lyapunov stable sets
\cite{bs} that
$\Lambda(Y)=\cap_n O_n$
where $O_n$ is a nested sequence of
positively invariant open
sets of $Y$. Obviously we can assume that
$O_n\subset U$ for all $n$.
Clearly the stable set of $O_n$ is
open in $U$.
Let us prove that such a stable set
is dense in $U$.
Let $O$ be an open subset
of $U$.
By Theorem \ref{th2} there is $x\in O$ such that
$\omega_Y(x)$ contains $\sigma(Y)$.
Clearly $\sigma(Y)$ belongs to
$O_n$
and so $\omega_Y(x)$ intersects $O_n$ as well.
Hence there is $t>0$ such that $X_t(x)\in
O_n$.
The last implies that $x$ belongs
to the stable set of $O_n$.
This proves that the stable set of $O_n$ is dense
for all $n$.
But the stable set of
$\Lambda(Y)$ is the intersection
of $W^s_Y(O_n)$ which is open-dense in $U$.
We conclude that the stable set
of $\Lambda(Y)$ is residual and the proof follows.
\end{proof}

Theorem \ref{tm} gives only a partial
answer for Question \ref{q1}
(in the one singularity case) since chain-transitive Lyapunov stable
set are not attractors in general.
However a positive answer for the question
will follow (in the one singularity case)
once we give positive answer for the
questions below.

\begin{q}
\label{q2}
Is a singular-hyperbolic, Lyapunov stable set an attracting set?
\end{q}

\begin{q}
\label{q3}
Is a singular-hyperbolic, chain-transitive,
attracting set a transitive set?
\end{q}

As it is well known these questions have positive answer
replacing singular-hyperbolic by hyperbolic
in their corresponding statements.
Besides it, a positive answer for
Question \ref{q2} holds provided the two
branches of the unstable manifold of every singularity
of the set are dense on the set \cite{mpa3}.

\medskip 

\flushleft
C. M. Carballo\\
Departamento de Matematica\\
Universidade Federal de Minas Gerais, ICEx - UFMG\\
Av. Antonio Carlos, 6627\\
Caixa Postal 702\\
Belo Horizonte, MG\\
30123-970\\
E-mail: carballo@mat.ufmg.br

\medskip

\flushleft
C. A. Morales\\
Instituto de Matem\'atica\\
Universidade Federal do Rio de Janeiro\\
P. O. Box 68530\\
21945--970 Rio de Janeiro, Brazil\\
E-mail: morales@impa.br

\end{document}